\documentclass[15pt,a4paper]{scrartcl}
\usepackage[utf8]{inputenc}
\usepackage[english,russian]{babel}
\usepackage{indentfirst}
\usepackage{misccorr}
\usepackage{graphicx}
\usepackage{amsmath}
\usepackage{amssymb}
\usepackage{graphicx}
\graphicspath{{pictures/}}
\DeclareGraphicsExtensions{.pdf,.png,.jpg}
\begin{document}

 \begin{center}
   {\Large Доказательство теоремы Радона при помощи понижения размерности}
\end{center}

 \begin{center}
   {\large Е.С.Колпаков}\footnote{Работа выполнена при частичной поддержке Добрушинской студенческой стипендии. 
Автор выражает благодарность А.Б. Скопенкову за полезные замечания и А. Акопяну за 
сообщение о статье
\cite{2}.}
\end{center}

Известна 
классическая теорема Радона [1]:

\smallskip

  \begin{it}Дано целое $d \geq 1$ и $d+2$ точек в $d$-мерном пространстве $\mathbb{R}^{d}$. Тогда эти точки можно разбить на два непересекающихся подмножества, чьи выпуклые оболочки имеют непустое пересечение.\end{it}

\smallskip

Первоначальное алгебраическое
доказательство  этой теоремы см. \cite{3}. Оно обычно и приводится. 
В настоящей статье 
содержится 
другое её
доказательство --- через понижение размерности. 
Ещё одно
доказательство через понижение размерности, 
несколько
более сложное, приводится в статье \cite{2}. 
В нём используется

\smallskip

 \textbf{Лемма об отделении}: Даны $d+2$ точек в $d$-мерном пространстве. Тогда существует гиперплоскость, натянутая на $d$ точек из данных, которая разделяет две оставшиеся точки.

\smallskip

 В настоящей статье лемма об отделении не используется. Оба доказательства через понижение размерности дают следующий, более сильный 

 Разбиение множества точек на два подмножества, выпуклые оболочки которых пересекаются {\it ровно по одной точке}, называется \textit{радоновским.} 

  \textbf{Количественная теорема Радона}. \begin{it} Для любого целого $d \geq 1$ и $d+2$ точек общего положения в $d$-мерном пространстве $\mathbb{R}^{d}$
радоновское
разбиение множества этих точек существует и единственно. \end{it}

 Теорема Радона связана со следующими теоремами, которые могут быть доказаны через понижение 
 размерности (см. \cite{4}, \cite{5}).

   \textbf{Теорема Конвея-Гордона-Закса}.  \begin{it} Для любых 
   шести
   точек в пространстве, никакие 4 из которых не лежат в одной плоскости, найдутся два зацеплённых треугольника с вершинами в этих точках. \end{it}

\smallskip

 \textbf{ Теорема ван Кампена-Флореса}.\begin{it} Среди любых 
 семи
 точек в четырёхмерном пространстве $\mathbb{R}^{4}$  можно выбрать две непересекающиеся тройки точек так, что треугольники с вершинами в них пересекаются. \end{it}

\smallskip
\smallskip

 {\it \textbf{Доказательство количественной теоремы Радона}}. 
Проведём индукцию
по $d$. 

{\it База:} $d=1$. Даны 3 точки на прямой. 
Тогда есть
ровно одно радоновское разбиение:  одно из множеств состоит из двух крайних точек, другое из средней точки.

 {\it Переход от $d-1$ к $d$. 
 } Обозначим
 через $A_1, ... , A_{d+2}$ данные точки в пространстве $\mathbb{R}^{d}$. Меняя, если нужно, нумерацию точек, можно выбрать гиперплоскость $\alpha$ так, чтобы точка $A_{d+2}$ лежала в одном полупространстве относительно $\alpha$, а точки $A_1, ... , A_{d+1}$ - в другом полупространстве. Тогда точка $A_{d+2}$ не лежит в симплексе $A_1...A_{d+1}$. Для $1 \le i \le d+1$ обозначим через $A_i^\prime$ пересечение $\alpha$ с $A_{d+2}A_i$. 

\begin{figure}[h!]
\includegraphics[scale=0.28]{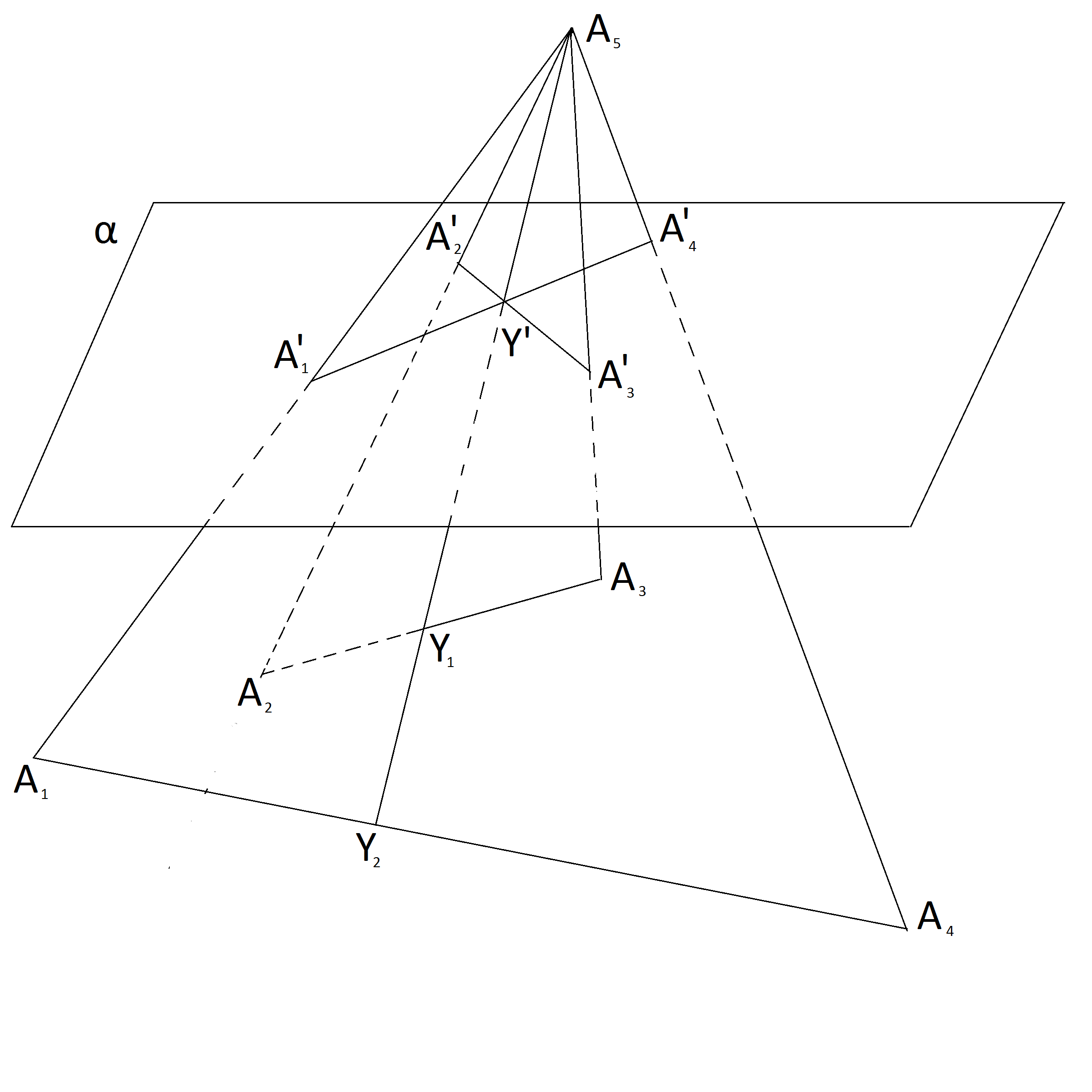}
\caption{ На рисунке изображён случай размерности $d=3$. Точка $A_5$ отделена плоскостью $\alpha$ от точек $A_1, A_2, A_3, A_4$. Отрезки $A_5A_1$, $A_5A_2$, $A_5A_3$, $A_5A_4$ пересекают плоскость $\alpha$ в точках $A_1^\prime, A_2^\prime, A_3^\prime, A_4^\prime$ соответственно. Отрезки $A_1^\prime A_4^\prime$ и $A_2^\prime A_3^\prime$ пересекаются в точке $Y'$. Луч $A_5Y'$ пересекает отрезки $A_2A_3$ и $A_1A_4$ в точках $Y_1$ и $Y_2$ соответственно.
Точка $Y_1$ лежит на отрезке $A_5Y_2$, поэтому отрезок $A_2A_3$ пересекает треугольник $A_1A_4A_5$ в точке $Y_1$.}
\end{figure}

Выпуклая оболочка множества $X$ будет обозначаться $\langle X \rangle$.

 1) {\it Доказательство
 существования радоновского разбиения.}
 По предположению индукции имеется 
 радоновское
 разбиение множества точек ${A_1^\prime, ... , A_{d+1}^\prime} = I^\prime \sqcup J^\prime$ в $(d-1)$-мерном пространстве $\alpha$. 
Поскольку точки $A_1, ... , A_{d+1}$ находятся в общем положении, пересечение $I'$ c $J'$ состоит из единственной точки, которую обозначим $Y'$
(рис. 1).

Пусть
${A_1, ... , A_{d+1}} = I \sqcup J$ соответствующее разбиение множества точек $A_1, ... , A_{d+1}$. 
А именно, пусть $I$ --- множество 
таких точек $A_i$,
что $A_i^\prime$ лежит в $I^\prime$. Аналогично определим множество $J$. 
Тогда отрезок $A_{d+2}Y'$ лежит в симплексах $\langle \{A_{d+2}\}\cup J \rangle$ и $\langle \{A_{d+2}\}\cup I \rangle$. Следовательно, луч $A_{d+2}Y'$ пересекает каждый из симплексов $\langle I\rangle$ и $\langle J \rangle$. Обозначим $Y_1 := A_{d+2}Y \cap \langle I \rangle$ и $Y_2 := A_{d+2}Y \cap \langle J  \rangle$.

 Не ограничивая общность рассуждения, считаем, что $Y_1$ лежит на отрезке $A_{d+2}Y_2$. Тогда симплекс $\langle I \rangle$ пересекает симплекс $\langle \{A_{d+2}\}\cup J \rangle$ по точке $Y_1$. Поскольку точки $A_1, ... , A_{d+2}$ находятся в общем положении, 
 получаем, что
 $Y_1$ - единственная точка пересечения $\langle I \rangle$ и $\langle \{A_{d+2}\}\cup J \rangle$.

Значит, $I \sqcup ( \{A_{d+2}\}\cup J)$ или $( \{A_{d+2}\}\cup I) \sqcup J$  является радоновским разбиением множества точек $A_1, ... , A_{d+2}$.

2. {\it Доказательство
единственности радоновского разбиения.}
 Предположим, что радоновских разбиений хотя бы два. 
 Пусть
 $I$ и $J$
 --- одно из них,
 $\widetilde{I}$ и $\widetilde{J}$                                                                                                                                                                                                                                                                                         
 --- другое. Положим
 $X=\langle I\rangle \cap \langle J \rangle$.

Без ограничения общности  
$A_{d+2} \in J$. 
Положим
$X^\prime = A_{d+2}X \cap \alpha$.  
Пусть
$I^\prime$ --- множество 
таких точек $A_i^\prime$,
что $A_i$ лежит в $I$, 
а
$J^\prime$ --- множество 
таких точек $A_i^\prime$,
что $A_i$ лежит в $J$. Тогда $X^\prime$ является единственной точкой пересечения  $\langle I^\prime \rangle$ и $\langle J^\prime \rangle$.  Значит, $I^\prime \sqcup J^\prime$ есть радоновское разбиение множества точек  $A_1^\prime, ... , A_{d+1}^\prime$.

 Аналогично из множеств $\widetilde{I}$ и $\widetilde{J}$ строится другое радоновское разбиение множества точек  $A_1^\prime, ... , A_{d+1}^\prime$, отличное от $I^\prime \sqcup J^\prime$.

 Таким образом, у множества точек $A_1^\prime, ... , A_{d+1}^\prime$ 
 из $(d-1)$-мерного пространства $\alpha$
 есть два радоновских разбиения, что противоречит предположению индукции. 
 Значит, число радоновских разбиений множества точек $A_1, ... , A_{d+2}$ не больше 1.
 Переход индукции доказан.

 $\linebreak$


\bigskip

\noindent Егор Сергееевич Колпаков, МГУ и НМУ \\
kolpak317bel@gmail.com

\end{document}